\documentclass[11pt]{article}
\usepackage{amsmath,amsthm,bbm,}

\newcommand{\R}{\mathbbm{R}}

\newcommand{\nothing}[1]{}

\newcommand{\De}{{\Delta}}
\newcommand{\na}{{\nabla}}
\newcommand{\pa}{{\partial}}

\newcommand{\la}{{\lambda}}

\newcommand{\om}{{\omega}}
\newcommand{\Om}{{\Omega}}

\newcommand{\g}{{\gamma}}

\newcommand{\al}{{\alpha}}
\newcommand{\bt}{{\beta}}
\newcommand{\te}{{\theta}}

\newcommand{\ep}{{\epsilon}}

\newcommand{\hr}{{\widehat{r}}}
\newcommand{\be}{{\begin{equation}}}
\newcommand{\ee}{{\end{equation}}}

\begin{document}
\title{\bf Open problems concerning the H\H older continuity of the direction of
vorticity for the Navier-Stokes equations. }
\author{H.~Beir{\~a}o~da Veiga,\\
Department of Mathematics,\\
Pisa University, Italy.\\
email: bveiga@dma.unipi.it}

\maketitle

\vspace{0.2cm}

\begin{abstract}
{In these notes we consider solutions $\,u\,$ to the evolution three
dimensional Navier-Stokes equations in the whole space. We set
$\om=\nabla\times u\,,$ the vorticity of $\,u\,.$ For brevity, we
set $\,L^s(\,L^r\,)=\,L^s(0,\,T;\,L^r(\R^3)\,)\,,$ and similar. Our
study concerns mainly the relation between $\bt-$H\H older
continuity assumptions on the direction of the vorticity and induced
regularity results like $\,\om \in\,L^{\infty}(\,L^r\,)\,.$ It is
well known that for $\,\bt=\,\frac12\,$ one gets $\,\om
\in\,L^{\infty}(\,L^2\,)\,.$ In particular, $\,u
\in\,L^{\infty}(\,L^6\,)\,$ follows. This shows that $\frac1{2}-$H\H
older continuity implies strong regularity for the velocity $\,u\,,$
by a classical result. On the other hand, it looks quite predictable
that a strictly decreasing perturbation of $\,\bt\,$ near $\,\frac12
\,,$ should induce a strictly decreasing perturbation for $\,r\,$
near $\,2\,.$ This possibility would imply strong regularity for
values $\,\bt <\,\frac12 \,,$ since strong regularity holds if
merely $\,\om \in \,L^{\infty}(\,L^r\,)\,,$ for some $\,r\geq
\,\frac32 \,.$ This result would be in contrast with a previous
conjecture which suggests that $\,\bt=\frac12\,$ is the smallest
value enjoying (in some non rigorous sense, see below) the above
strong regularization property. In fact, in the sequel, the smallest
value $\,\bt\,$ for which we will be able to prove that $\,\om
\in\,L^s(L^r\,)\,$ for any $\,r \in (1,\,2]\,,$ is still
$\,\bt=\,\frac12\,.$ This conclusion reenforces the above conjecture
on the particular significance of the value $\,\bt=\,\frac12\,.$
But, on the other hand, it also shows the lack of a complete
understanding of the full phenomena since the above perturbation
argument looks well-founded. In the last section we discuss some
related open problems, consistent with our calculations. We hope
that the approach bellow should help readers interesting in carrying
on this investigation.}
\end{abstract}
\vspace{0.2cm}


\vspace{0.2cm}

\section{Introduction. Motivation and conclusions.}\label{uns}
In the following we consider the evolution Navier-Stokes equations
\begin{equation}\label{nse}
\left\{  \begin{aligned} u_t+(u\cdot\nabla)\,u-\,\triangle u+\nabla
p&=0\,,\\
\nabla\cdot u&=0\qquad\text{in } \R^3 \times (0,T], \\
u(x,0)&=u_0(x)\qquad\text{in } \R^3\,,
\end{aligned}\right.
\end{equation}
where the divergence free initial data $\,u_0\,$ belongs (for
instance) to the Sobolev space $\,H^1(\R^3)\,.$ For simplicity, we
assume that external force vanishes, and viscosity is equal to one.
If we want solutions which are regular including the time
$\,t=\,0\,,$ then $\,u_0\,$ should be assumed more regular. This
kind of problem is out of real interest here. We consider the
simplest situation, the whole space case, even if no-slip boundary
conditions could be considered.\par%
We will not repeat well know notation as, for instance, Sobolev
spaces notation, and so on. For brevity, we set
$\,L^s(\,L^r\,)=\,L^s(0,\,T;\,L^r(\R^3)\,)\,,$ and similar.
Solutions $\,u\in\,L^2(0,\,T;\,H^1(\R^3)\,)
\cap\,L^{\infty}(0,\,T;\,L^2(\R^3)\,)\,$ are defined in the well
known Leray-Hopf weak sense. In the sequel, for brevity, we say that
a weak solution $u$ is \emph{strong} if
\begin{equation}
\label{queri}%
u \in L^s (0,\,T; L^q(\Om)\,)\,,
\end{equation}
for some exponents $s$ and $q$, $\,q\geq\,3\,,$ satisfying
$\,\la(s,\,q) \leq\,1\,,$ where
\begin{equation}
\label{quereum}%
\la(s,\,q) \equiv \,\frac2{s}+\,\frac3{q}\,.
\end{equation}
In the following $\omega=\nabla\times u\,$ denotes the
\emph{vorticity} field. By a well known Sobolev's embedding theorem
equation \eqref{queri} and equation
\begin{equation}
\label{queriom}%
\om \in L^s (0,\,T; L^r(\Om)\,)
\end{equation}
are related by
\begin{equation}
\label{qqrr}%
\frac1{q}=\,\frac1{r}-\,\frac13 \,.
\end{equation}
Consequently, solutions are strong if \eqref{queriom} holds for
\begin{equation}
\label{pich}%
\frac2{s} +\,\frac3{r} \leq\,2\,.
\end{equation}
In particular, solutions are strong if
\begin{equation}
\label{elrr}%
\om\in\,L^{\infty}(0,\,T;\,L^r(\Om)\,)%
\end{equation}
for some $\,r \geq\,\frac32\,.$ This point has a main role in the
sequel.

\vspace{0.2cm}

Let's go inside the core of this notes. We set
$$
\theta(x,y,t)\overset{def}{=}\angle(\omega(x,t),\,\omega(y,t))\,,
$$
where the symbol "$\,\angle\,$" denotes the amplitude of the angle
between two vectors. We are interested in studying possible
regularization effects of $\bt-$H\H older continuity assumptions on
the direction of the vorticity, namely
\begin{equation}\label{condizione 2}%
\sin\theta(x,y,t)\leq\,\, c\,|x-y|^\beta
\end{equation}
in $\,\R^3\times\,\R^3\times(0,T)\,,$ for some
$\,\beta\in(0,1/2]\,.$\par%
Our study concerns mainly the possible relations between the above
assumption and regularity results like \eqref{elrr}. It is well
known, see \cite{bb2002}, that assumption \eqref{condizione 2} with
$\,\bt=\,\frac12\,$ implies \eqref{elrr} with $\,r=\,2\,$ (this
situation is called here the \emph{Hilbertian case}). Strong
regularity follows. However, as already remarked, strong regularity
is still guaranteed by \eqref{elrr}, if merely $\,r\geq \,\frac32
\,.$ This shows the weight of being able to prove \eqref{elrr} for
some $\,r\geq \,\frac32 \,.$ It is worth noting that this result
looks quite natural. In fact, a strictly decreasing perturbation of
$\,\bt\,$ near $\,\frac12\,$ should induce a corresponding strictly
decreasing perturbation of $\,r\,$ near $\,2\,.$ In particular,
values $\,r \in (\frac32,\,2)\,$ would be obtained. Consequently
strong regularity would follow for $\,\bt<\,\frac12\,.$ But, on the
opposite side, this regularity result would be in some contrast with
the opinion supported by us in a previous publication, see the
appendix in reference \cite{bdv2014}, which predicts that strong
regularity under the $\,\bt=\,\frac12\,$ assumption should be "as
strong as" the classical sufficient condition \eqref{queri} for
$\,\lambda=\,1\,.$ Roughly speaking, to improve the
$\,\bt=\,\frac12\,$ sufficient condition for strong regularity could
be as hard as to extend the sufficient condition \eqref{queri} to values $\,\lambda>\,1\,.$%

\vspace{0.2cm}

In the sequel we follow the perturbation strategy refereed above,
namely, generating a strictly decreasing perturbation of $\,\bt\,,$
in the proof of the Hilbertian case
$\,(\bt,\,r)=\,(\,\frac12,\,2)\,,$ and waiting for a corresponding
drop of the value of the exponent $\,r\,.$ This strategy requires to
have in hands a suitable extension of the known $\,r=\,2\,$ proof,
which works for values $\,r \neq\,2\,.$ Suitable means here that the
case $r=\,2\,$ must be "perfectly" embedded in the new proof.
Roughly speaking, with "continuity of argumentation" with respect to
the variable $\,r\,.$ Actually the extension below does not require
$\,r \neq\,2\,.$ The construction of this more general proof is a
main point in the sequel.\par%
It is worth noting that the smallest value $\,\bt\,$ obtained by us
below, which guarantees \eqref{elrr} regularity for exponents $\,r
\in (1,\,2]\,,$ is still $\,\bt=\,\frac12 \,.$ We believe that this
"negative" conclusion is not due to an insufficiency of the above
extended new proof, but is due to more subtle reasons. From one
side, the conclusion reenforces our previous claim in favor of the
real, not merely technical, significance of the particular value
$\,\bt=\frac12\,$. On the other side, it also shows a partial
failure of a very natural attempt to improve the $\,\bt=\frac12\,$
sufficient condition for regularity, by means of a perturbation
approach. Summing up, a good understanding of the full phenomena
persists. In particular, the link between $\,\bt-$H\H older
continuity assumptions on the direction of the vorticity, for $\,\bt
<\,\frac12 \,,$ and regularity results of type
$\,\om\in\,L^s(\,L^r\,)\,$ is still an open problem. We hope that
the specific approach presented bellow should help readers
interested in carrying on this challenging investigation.
\section{Some related known results.}\label{unspoc}
We start by some references, strongly related to the present notes
by methods of proof. We begin by recalling the very fundamental
pioneering paper \cite {CF93}, by P.Constantin and Ch.Fefferman,
where the authors prove, in particular, that solutions to the
evolution Navier-Stokes equations in the whole space are smooth if
the direction of the vorticity is Lipschitz continuous with respect
to the space variables, namely assumption \eqref{condizione 2} for
$\,\bt=\,1$. This condition is assumed for almost all $\,x$ and $y$
in $\,R^3\,,$ and almost all $t$ in $(0,\,T)\,.$ Actually, in \cite
{CF93}, the assumption is merely required for points $x$ and $y$
where the vorticity at both $x$ and $y$ is larger than a given,
arbitrary constant $\,k.$ This improvement was, or can be, extended
in the same way to many subsequent papers on the subject. It is also
easily show that assumption \eqref{condizione 2} can be restricted
to couples of points $x$ and $y$ satisfying $|x-y|<\,\delta,\,$ for
an arbitrary positive constant $\delta$.\par%
In reference \cite{bb2002}  L.C. Berselli and the present author
showed, in particular, that regularity still holds in the whole
space by replacing Lipschitz continuity by $\,\frac12-$H\H older
continuity. This is, up to now, the strongest result in the
literature. In the subsequent paper \cite{olga}, by a trivial
modification of the proof shown in \cite{bb2002}, the following
result is proved. Let $\,\beta\in(0,1/2]\,$ be given. Further,
assume that condition \eqref{condizione 2} holds and, in addition,
that $\,\omega \in L^2(L^r)\,$ where $\,r=\frac{3}{\beta + 1}\,$.
Under these assumptions it was shown that the solution $u$ is
strong. In particular, the sufficiency of condition
$\,\bt=\,\frac12\,$ for regularity is reobtained, since $\,\omega
\in L^2(L^2)\,$ holds for any weak solution. This result is
particularly related to the problem treated in the sequel.

\vspace{0.2cm}

Concerning other related papers, we start by recalling reference
\cite{bdvCPAA} where we extended the above kind of results to the
half-space $\,\Om=\,\R^3_+\,$, endowed with the slip boundary
condition (``stress-free'' boundary condition)
\begin{equation}
\label{bc}%
\left\{  \begin{aligned}
u\cdot n&=0\,, \\
\omega\times n&=0\,,\\
\end{aligned}\right.
\end{equation}

where $n$ is normal to the boundary. In reference \cite{bvbers},
L.C. Berselli and the author succeed in extending this result to the
case in which $\Omega\subset \R^3\,$ is an open, bounded set with a
smooth boundary, by appealing to suitable representation formulas
for Green's matrices. In reference \cite{luigi} regularity is proved
by replacing continuity requirements on $\,\sin\theta(x,y,t)\,$ by a
smallness assumption. Essentially, it is proved that there is a
sufficiently small constant $\,C_1\,$ such that regularity holds if
$\, \sin\theta(x,y,t)\leq \,C_1\,.$
Clearly, there are many very interesting papers related to our
contributions. We recall here, without any claim of completeness,
references \cite{bdvmagenes}, \cite{bdv2014}, \cite{MR2563656},
\cite{cordoba}, \cite{MR2351138}, \cite{MR2645147},
\cite{MR2338368}, \cite{MR2436725}, \cite{CFM}, \cite{MR1972658},
\cite{MR3010194}, \cite{MR3026558}, \cite{MR2782615},
\cite{MR2525642}, \cite{MR2669442}, \cite{MR2095449},
\cite{MR2202302}, \cite{MR2236569} \cite{MR2232429},
\cite{MR2062645}, \cite{vasseur}.
\section{The main new estimate.}\label{tres}%
In the following $\,f(s)\,$ denotes a real continuous differentiable
function $\,f:\,\R^+ \rightarrow \R^+\,.$ We set
$$
F(s)=\, \int_0^s \,f(\tau) \,d\tau\,.
$$
Hence $\,F'(s)=\,f(s)\,.$ By applying the curl operator to equation
\eqref{nse} we get the
well-known equation%
$$
\omega_t+(u\cdot\nabla)\,\omega-\,\nu\,\Delta \omega
=(\omega\cdot\nabla)\,u\,.
$$
Scalar multiplication by $\,f(|\om|^2)\,\om\,$, integration in
$\,\R^3\,$, and integrations by parts easily show that
\begin{equation}\label{enerotum}
\frac12\,\frac{d}{dt}\,\int\,F(|\om|^2)\,dx -\int\,f(|\om|^2)\,\De
\,\om \cdot\,\om \,dx =\,\int\,f(|\om|^2)\,(\omega\cdot\nabla)\,u
\cdot\,\om \,dx\,.
\end{equation}
Non-labeled integrals are over $\,\R^3\,$. Straightforward
calculations show that
$$
-\int\,f(|\om|^2)\,\De \,\om \cdot\,\om \,dx =\int\,f(|\om|^2)\,|\na
\,\om|^2 \,dx +\,2\,\int\,f'(|\om|^2)\, |\om|^2\,|\na \,\om|^2
\,dx\,.
$$
Hence
\begin{equation}\label{enerotum}
\begin{split}
\frac12\,\frac{d}{dt}\,\int\,F(|\om|^2)\,dx
+\,\int\,f(|\om|^2)\,|\na \,\om|^2 \,dx \phantom{aaaaaaaaaaaaaa}\\
\\
\leq\,2\,\int\,f'(|\om|^2)\, |\om|^2\,|\na \,\om|^2 \,dx
+\,\int\,f(|\om|^2)\,(\omega\cdot\nabla)\,u \cdot\,\om \,dx\,.
\end{split}
\end{equation}

\vspace{0.2cm}

In these notes we are interested in the particular case
$\,f(s)=\,s^{-\,\al}\,.$ In the Hilbertian case, in which
$\,\al=\,0\,,$ many of the devices used in the sequel are superfluous.\par%
It is useful to start by considering the approximation functions
\begin{equation}%
\label{effepsilon}%
f_{\ep}(s)=\,(\ep+\,s)^{-\,\al}\,,
\end{equation}
where $\,\ep>\,0\,,$ and $\,0 \leq \,\al \leq \,\frac12 \,.$ Note
that $\,f'_{\ep}(s)<\,0\,.$ Straightforward calculations show that
the absolute value of the first integral on the right hand side of
equation \eqref{enerotum} is bounded by $\,\al\,$ times the second
integral on the left hand side of the same equation. So, one has

\begin{equation}
\label{enerotdois}
\begin{split}
\frac{1}{2(1-\al)}\,\frac{d}{dt}\,\int\, (\ep+\,|\om|^2)^{1-\,\al}
\,dx +\,(1-2\,\al) \int\,(\ep+\,|\om|^2)^{-\,\al} |\na \,\om|^2
\,dx\\
\\
\leq\, \int\, (\ep+\,|\om|^2)^{-\,\al}\,|\mathcal{K}(x)|\,dx%
\end{split}
\end{equation}
where
\begin{equation}
\label{kapa}%
\mathcal{K}(x):= \,\left(\,( \omega\cdot\nabla)\,u \cdot\,\omega
\right)(x)\,\,.
\end{equation}
Next we estimate from below the second integral on the left hand
side of equation \eqref{enerotdois} (see \cite{bdv1995} and
\cite{bdv87} for similar manipulations).  One has
\begin{equation}
\label{commod}%
(\ep+\,|\om|^2)^{-\,\al}\, |\na \,|\om|\,|^2\,
=\,\frac{1}{(1-\al)^2}\,\frac{|\om|^{2\,\al}}{(\ep+\,|\om|^2)^\al\,
}\, \big|\,\na(\,|\om|^{1-\,\al})\,\big|^2\,.
\end{equation}
Since $|\,\na\,\om|\,\geq |\,\na\,|\om|\,|\,,$ it follows from
equation \eqref{enerotdois} that
$$
\begin{array}{l}
\frac{1}{2(1-\al)}\,\frac{d}{dt}\,\int\,(\ep+\,|\om|^2)^{1-\,\al}
\,dx +\,\frac{(1-2\,\al)}{(1-\al)^2}\, \int\, \big|\,\na(\,|\om|^{1-\,\al})\,\big|^2 \,dx\\
\\
\leq\, \int\, (\ep+\,|\om|^2)^{-\,\al}\,|\mathcal{K}(x)|\,dx\,.%
\end{array}
$$
By letting $\,\ep \rightarrow \,0$ one gets
\begin{equation}
\label{enerotress}%
\begin{array}{l}
\frac{1}{2(1-\al)}\,\frac{d}{dt}\,\|\om\|^{2(1-\al)}_{2(1-\al)}
+\frac{(1-2\,\al)}{(1-\al)^2}\,\|\na(\,|\om|^{1-\,\al})\|^{2}_2 \\
\\
\leq\, \int\, |\om|^{-2\,\al}\,|\mathcal{K}(x)|\,dx%
\end{array}
\end{equation}
which, for $\,\al=\,0\,,$ is precisely the estimate obtained in the
Hilbertian case. Now we recall that, in the Hilbertian case, one
appeals to the Sobolev's embedding $\, \|g\|_6
\leq\,c_0\,\|\na\,g\|_2\,$ to allow substitution of
$\,\|\na\,\om\|_2\,$ by $\,\|\om \|_6\,.$ In the more general case
considered here we appeal to the same device, by  applying the above
Sobolev's estimate to the function $\,g=\,|\om|^{1-\,\al}\,.$  After
this device, equation \eqref{enerotress} reads
\begin{equation}
\label{enerotres}%
\frac{d}{dt}\,\|\om\|^{2(1-\al)}_{2(1-\al)}
+\,c_1\,\|\om\|^{2(1-\al)}_{6(1-\al)}
\leq\, \int\, |\om|^{-2\,\al}\,|\mathcal{K}(x)|\,dx\,.%
\end{equation}
The symbol $c\,$, and similar, may denote distinct positive
constants.
\section{Estimating the nonlinear term by a Riesz potential.}\label{tresete}%
In this section we estimate the left hand side of equation
\eqref{enerotress} by means of a related Riesz potential. This is
one of the main ideas introduced by Constantin and Fefferman in
\cite{CF93}, and took again in \cite{bb2002}. We follow here the
presentation given in reference \cite{bdv2014} (where bounded
domains are considered). See also \cite{bvbers}.\par%
Since $\, -\,\Delta\,u=\,\nabla \times\,(\nabla \times u)-\, \nabla
\,(\nabla\cdot\,u)\,,$ it follows that
\begin{equation}\label{lapu}
-\Delta\,u=\,\nabla \times \omega \quad \textrm{in} \quad \R^3\,,
\end{equation}
for each $t$. So
\begin{equation}
\label{ugreen}%
u(x)=\, \int \,G(x,y)\,(\nabla \times \omega)(y) \,dy\,,
\end{equation}
where
$$
G(x,y)=\,\frac{1}{4\,\pi\,|x-\,y|}\,.
$$
In particular
\begin{equation}
\label{base}%
\left|\frac{{\partial}^2\,G(x,y)}{\partial\,y_k\partial\,x_i}\right|\leq\,\frac{c}{|x-\,y|^3}\,.
\end{equation}
 Set, for each triad $(j,k,l)$, $j,k,l \in
\,\left\{1,2,3\right\}$,
$$
\epsilon_{ijk}= \left\{ \begin{array}{ll}%
1  & \textrm{if $(i,j,k)$ is an even permutation}\,,\\
-1 & \textrm{if $(i,j,k)$ is an odd permutation}\,,\\
0 & \textrm{if two indexes are equal}\,.%
\end{array}%
\right.%
$$
These are the components of the totally anti-symmetric Ricci tensor.
One has
\begin{equation}\label{psterno}
(a \times b)_j=\,\epsilon_{jkl}\,a_k\,b_l\,, \qquad (\nabla \times
v)_j=\,\epsilon_{jkl}\,\partial_k v_l\,.
\end{equation}
The usual convention about summation of repeated indexes is assumed.\par%
In particular
\begin{equation}
\label{base}%
\left|\frac{{\partial}^2\,G(x,y)}{\partial\,y_k\partial\,x_i}\right|\leq\,\frac{c}{|x-\,y|^3}\,.
\end{equation}
By considering in equation \eqref{ugreen} a single component
$u_j\,,$ and by appealing to \eqref{psterno}, an integration by
parts yields
$$
u_j(x)=\,\int \,G(x,y)\,\epsilon_{jkl}\,\pa_k\, \omega_l(y)\,dy =\,
-\int \epsilon_{jkl}\, \frac{\partial\,G(x,y)}{\partial\,y_k}\,
\omega_l(y)\,dy\,.
$$
Hence
$$
\frac{\partial\,u_j(x)}{\partial\,x_i}= -P.V. \,\int
\epsilon_{jkl}\,\frac{{\partial}^2\,G(x,y)}{\partial\,x_i\partial\,y_k}\,
\omega_l(y)\,dy \,.
$$
It readily follows that
$$
\begin{array}{l}
\mathcal{K}(x) =\,-\int
\epsilon_{jkl}\,\frac{{\partial}^2\,G(x,y)}{\partial\,y_k\partial\,x_i}\,
\omega_i(x)\,\omega_j(x)\, \omega_l(y)\,dy\,.
\end{array}
$$
Since $\,-\,\epsilon_{jkl}\,\omega_j(x)\,\omega_l(y)=\,
\big(\,\omega_j(x) \times \,\omega_l(y)\,\big)_k\,,$ it follows that
$$
\mathcal{K}(x)= \,P.V. \,\int \,
\frac{{\partial}^2\,G(x,y)}{\partial\,y_k\,\partial\,x_i}\,
\omega_i(x)\,\big(\,\omega_j(x) \times \,\omega_l(y)\,\big)_k\,
\,dy\,.
$$
By appealing to \eqref{base} one shows that
\begin{equation}\label{eeu}%
|\,\mathcal{K}(x)\,| \leq \, \int \,\frac{c}{|x-\,y|^3}\,
|\omega(x)|^2 \,|\omega(y)|\,\sin \theta(x,y,t) \,dy \,.
\end{equation}
Now we appeal to the main assumption \eqref{condizione 2}, where
$\,\beta\in[0,1/2]\,.$ By appealing to \eqref{eeu} one gets
\begin{equation}
\label{mmaj3}%
|\mathcal{K}(x)|\leq \,c\,|\omega(x)|^2\, I(x)\,,
\end{equation}
where
$$
I(x)=\,\int_{\Omega} \, |\omega(y)|\,\frac{dy}{|x-\,y|^{3-\,\beta}}
$$
is the Riesz potential in $\,\R^3\,.$ Recall that (see
\cite{MR44:7280}) if $\,0<\,\beta<\,3\,,$ and if $\,\omega \in
\,L^{\hr}(\Omega)\,$ for some $\,1<\,\hr<\,3 \,,$ then
\begin{equation}\label{tomas}%
\|I\|_q \leq \,c\,\|\omega\|_{\hr}\,,
\end{equation}
where
$$
1/q=\, 1/\hr \,-\,\beta/\,3\,.
$$
In particular, by \eqref{mmaj3}, the right hand side of equation \eqref{enerotres} satisfies the estimate%
\begin{equation}
\int\, |\om|^{-2\,\al}\,\mathcal{K}(x)\,dx \leq\,c\,\int
|\om|^r\,I(x)\,dx\,,
\end{equation}
where $\,r=\,2(1-\,\al)\,.$ So,
\begin{equation}
\label{abas}%
\frac{d}{dt}\|\omega\|_r^r+\,\|\omega\|^r_{3r} \leq\,c\, \int
|\om|^r\,I(x)\,dx\,,
\end{equation}
for
$$
1<r\leq\,2\,.
$$
From now on we eliminate the above parameter $\al$ by appealing to
the new exponent $r\,.$ Next, by appealing to \eqref{tomas}, we
write
\begin{equation}\label{abas2}%
\frac{d}{dt}\|\omega\|_r^r+\,\|\omega\|^r_{3r}
\leq\,c\,\|\om\|_{\hr}\, \|\om\|^{r}_{q'r}\,\,,
\end{equation}
where
\begin{equation}
\label{qlinha}%
\frac1{q'}=\, 1-\,\frac1{\hr} \,+\,\frac{\beta}3\,.
\end{equation}
More precisely, by \eqref{enerotress}, we could write (not used in
the sequel)
\begin{equation}\label{abastat}%
\frac{d}{dt}\|\omega\|_r^r+\, \|\,\na\,|\omega|^{\frac{r}2}\,\|^2_2
\leq\,c\,\|\om\|_{\hr}\, \|\om\|^{r}_{q'r}\,\,.
\end{equation}
\section{Conclusions.}\label{conclusions}%
Our final task, which is the main task in this section, is looking
for pairs $r$ and $\bt$ such that the typical estimate
\begin{equation}\label{tens}%
\|\om\|_{\hr}\,\|\om\|^{r}_{q'r} \leq\, C_\ep
\,\|\om\|^2_2\,\|\om\|^r_r +\,\ep\, \|\,\om\|^r_{3r}
\end{equation}
holds. As usual, the meaning of the above equation is that $\ep\,$
may be any positive, arbitrarily small real number, at the price of
having a corresponding, arbitrarily large, value of $C_\ep\,.$ The
motivation for this requirement is standard. Assume that
\eqref{tens} holds for some pair of values $r$ and $\bt$. Then, by
appealing to \eqref{abas2}, to \eqref{tens}, to
\begin{equation}
\label{crucia}%
\|\omega(t)\|^2_2 \in\,L^1(0,\,T)\,,
\end{equation}
and to Gronwall's lemma, we show that
\begin{equation}
\label{infasr}%
\om \in L^\infty (0,\,T; L^r(\Om)\,)\cap\, L^r (0,\,T;
L^{\,3r}(\Om)\,)\,.
\end{equation}
Even though $r \neq 2\,,$ a central role in the right hand side of
equation \eqref{tens} is still required to the integrability
exponent $2\,$. The reason for this choice is that \eqref{crucia} is
the strongest known estimate for the vorticity of weak solutions.\par%
Our aim is now to find pairs $\,(r,\,\bt)\, \in
\,(1,\,2]\times\,(0,\,1/2 \,]\,$ such that \eqref{tens} holds. It is
convenient to start by recalling the way followed in the proofs in
the Hilbertian case $\,r=\,2\,.$ In this case one had
$\,(r,\,\bt,\,\hr\,,q,\,q')=\,(2,\,1/2,\,2,\,3,\,3/2\,)\,.$ Hence
equation \eqref{abas2} reads
\begin{equation}
\label{fucas}%
\frac{d}{dt}\|\omega\|_2^2+\,\|\omega\|^2_6 \leq\,c\,\|\om\|_2\,
\|\om\|^{2}_3\,.%
\end{equation}
The next move in the classical proof was to appeal to the
\emph{interpolation} inequality
$$
\|\,\om\,\|^2_3 \leq\,\|\,\om\,\|_2 \,\|\,\om\,\|_6
$$
which, together with \eqref{fucas}, leads to the desired estimate
$$
\frac{d}{dt}\|\omega\|^2_2 +\,\|\omega\|^2_6 \leq\,C_\ep
\,\|\om\|^2_2\,\|\om\|^2_2\,+\, \ep\, \|\om\|^{2}_6\,.%
$$
Following the same idea, we decompose both norms $\,\|\om\|_{\hr}\,$
and $\,\|\om\|^{r}_{q'r}\,$ in the right hand side of \eqref{abas2},
by appealing to interpolation. Note that, in the Hilbertian case,
the choice $\hr=\,2\,$ in \eqref{tens} looks obvious. On the
contrary, in the more general situation treated here, a previous
restriction could cut potential to the proofs. So we opted here for
giving the largest possible width to the range of the parameter
$\,\hr\,$ by assuming that
$$
r\leq \,\hr\,\leq\,3r\,,
$$
as suggested by the left hand side of \eqref{tomas}. The same
freedom will be given to $\,q'\,r\,.$ More precisely, we start by
considering parameters $\,\al,\,\te,\,\g\,$ and
$\,\al',\,\te',\,\g'\,,$ in the interval $\,[0,\,1]\,,$ satisfying
$\,\al+\te+\g=\,\al'+\te'+\g'=\,1\,,$ and related to the exponents
$\,q'\,r\,$ and $\,\hr\,$ by the following equations:
\begin{equation}
\label{follas}%
\left\{\begin{aligned}%
\frac1{q'r}=\,\frac{\al}{r}+\,\frac{\te}2+\,\frac{\g}{3r}\,,\\
\frac1{\hr}=\,\frac{\al'}{r}+\,\frac{\te'}2+\,\frac{\g'}{3r}\,.
\end{aligned}\right.
\end{equation}
It follows, by interpolation, that
\begin{equation}
\label{dec1}%
\left\{\begin{aligned}%
\|\om\|_{q'r} \leq\,\|\om\|^{\al}_r \,\|\om\|^{\te}_2\,
\|\om\|^{\g}_{3r}\,,\\
\|\om\|_{\hr} \leq\,\|\om\|^{\al'}_r \,\|\om\|^{\te'}_2\,
\|\om\|^{\g'}_{3r}\,.
\end{aligned}\right.
\end{equation}
The values of the above parameters will be fixed in the sequel. One
has
\begin{equation}
\label{filas}%
B=: \|\om\|^{r}_{q'r}\,\|\om\|_{\hr}\leq\,\|\om\|^{\al'+\al\,r}_r
\,\|\om\|^{\te'+\,\te\,r}_2\, \|\om\|^{\g'+\,\g\,r}_{3r}\,.
\end{equation}
Next, by appealing to the dual exponents
$$
\frac{r}{\g'+\,\g\,r}\,, \quad  \frac{r}{(1-\,\g) r -\g'}\,,
$$
we get
\begin{equation}
\label{files}%
B\leq\, C_\ep \, \|\om\|^{\frac{(\al'+\al\,r)\,r}{(1-\,\g)r
-\g'}}_r\, \|\om\|^{\frac{(\te'+\te\,r)\,r}{(1-\,\g)r -\g'}}_2
+\,\ep\, \|\om\|^{r}_{3r}\,.
\end{equation}
We want
\begin{equation}
\label{querer}%
\frac{\al'+\al\,r}{(1-\,\g)\,r-\,\g'} =\,1\,, \quad
\frac{\te'+\te\,r}{(1-\,\g)\,r-\,\g'} =\,\frac2r\,,
\end{equation}
since this leads to
\begin{equation}\label{devac}%
\frac{d}{dt}\|\omega\|_r^r+\,\|\omega\|^r_{3r} \leq\,C_\ep\,
\|\om\|^2_2\,\|\om\|^r_r+\,\ep\,\|\omega\|^r_{3r}\,.
\end{equation}
By setting $\,\g=\,1-(\al+\te)\,$ and $\,\g'=\,1-(\al'+\te')\,$ in
the first equation \eqref{querer}, one easily shows that
\eqref{querer} is equivalent to
\begin{equation}\label{nsen}
\left\{ \begin{aligned} \te'+\te\,r=\,1\,,\\
\al'+\al\,r=\,\frac{r}2\,.
\end{aligned}\right.
\end{equation}
In addition, the exponents $\,q'\, $ and $\,\hr\,$ must verify
equation \eqref{qlinha}. This requirement is easily satisfied since
the parameter  $\bt\,$is still free. In other words, it gives the
value of the H\H older exponent $\bt=\,\bt(r)\,$ that leads to the
regularity result \eqref{infasr}. Let's calculate this value. By
appealing to the equation \eqref{qlinha} and to the second equation
\eqref{follas}, one shows that the first equation \eqref{follas} can
be written in the equivalent form
$$
1-\,\frac{\al'}{r}-\,\frac{\te'}2-\,\frac{\g'}{3r}+\,\frac{\bt}{3}=\,
\al+\,\frac{r}{2}\,\te +\,\frac{\g}3\,.
$$
Further, by replacing $\g$ and $\g'$ respectively by
$\,1-\,(\al+\,\te)\,$ and $\,1-\,(\al'+\,\te')\,$, straightforward
calculations lead to the desired expression of $\,\bt(r)\,,$ namely
\begin{equation}
\label{poder}%
\bt(r)=\,\frac2r\, (\al'+\al\,r) +\, \big(\,\frac32 -\,\frac1r
\,\big) (\te'+\te\,r)-\,2+\,\frac1r\,.
\end{equation}
Lastly, by appealing to \eqref{nsen} and \eqref{poder}, we realize,
in agrement to our prediction but also with some disappointment,
that the exponent $\,\bt\,$ obtained here does not depend on
$\,r\,.$ Actually, one gets
\begin{equation}\label{semvac}%
\bt=\,\frac12\,.
\end{equation}
Summing up, our attempt to feel out if the regularity result
\eqref{infasr} may hold for $\,r<2\,$ under a $\bt-$H\H older
continuity assumption on the direction of the vorticity, for some
$\bt<\,\frac12\,,$ has had here a negative reply. This conclusion
supports the argument, defended in the appendix of \cite{bdv2014},
that values $\bt<\,\frac12\,$ "does not imply strong regularity".%

\vspace{0.2cm}

Note that \eqref{nsen} shows that "natural" pairs of values are
$\al'=\,0\,$, $\,\al=\,\frac12\,,$ and $\,\te'=\,1\,,$
$\,\te=\,0\,$. Consequently, $\,\g=\,\frac12\,,$ and
$\,\al'=\,\g'=\,0\,.$ These are exactly the values used in the proof
of the Hilbertian case $\,r=\,2\,$.
\section{An open problem.}\label{abertas}%
In these notes we have considered regularity results in the form
shown by equation \eqref{elrr}, hence with $\,s=\,\infty\,.$ It
could be that the desired additional regularity holds in the
framework of finite values of $\,s\,.$ A particular case would be to
show that
$$
\om \in L^2 (0,\,T; L^r(\Om)\,)\,,
$$
for some $\,r>\,2\,.$ Weak solutions satisfy this requirement for
$\,r=\,2\,.$ The following is a more general open problem:\par%
Assume that a $\bt-$H\H older continuity assumption on the direction
of the vorticity holds for some value $\bt<\,\frac12\,.$ To show
that \eqref{queri} holds for a couple of exponents $s$ and $q$
satisfying
\begin{equation}
\label{nicas}%
\la(s,\,q)=\,\frac2{s} +\,\frac3{q}<\,\frac32 \,.
\end{equation}%
Note that $\,\lambda=\,1\,$ means smoothness. The regularity
required by equation \eqref{queri}, if $\,\lambda(s,\,q)=\,\frac32
\,,$ corresponds to the maximum regularity known for generical weak
solutions. Hence \eqref{nicas} would show some additional regularity
to weak solutions, in terms of integrability, due to the above
regularity assumption in terms of direction of vorticity. The
possibility of this "transfer" of regularity looks quite natural, by
taking into account that it holds in the Hilbertian case. Arriving
to a reply to the above open problem, positive or negative, is a challenging open problem.%

\vspace{0.2cm}

Equations \eqref{queri} and \eqref{queriom} are related by
\eqref{qqrr}. Roughly, the assumption \eqref{queri} with exponents
satisfying \eqref{nicas}, is "near equivalent" to the assumption
\eqref{queriom} with exponents satisfying
\begin{equation}
\label{quereo}%
\frac2{s} +\,\frac3{r}<\,\frac52 \,.
\end{equation}
However, it could be not negligible that \eqref{queriom} is
stronger.

\vspace{0.2cm}

We may also consider reverse problems, like to verify whether
regularity assumptions of type \eqref{queri} imply, or not imply,
regularity for the direction of the vorticity. Our guess is that
there is not a "clean" reply to this problem.

\end{document}